\documentclass[12pt,a4paper]{article}
\usepackage[utf8]{inputenc}
\usepackage{lmodern}
\usepackage[T1]{fontenc}
\usepackage[english]{babel}
\usepackage{ifpdf}
\usepackage[usenames,dvipsnames]{xcolor}
\usepackage{graphicx}
\usepackage[shortlabels]{enumitem}
\usepackage[numbers,sort]{natbib}
\usepackage[font=footnotesize,width=\textwidth]{caption}
\usepackage{booktabs}
\usepackage{hyperref}
\usepackage{uri}
\usepackage{tikz}
\usepackage{float}
\usepackage[top=1.8cm, bottom=2.8cm]{geometry}
\usepackage[parfill]{parskip}
\usepackage[affil-it]{authblk}
\usepackage[setbuildercolon]{matmatix}

\newcommand{\T}{\mathcal T}
\newcommand{\EE}{\mathbb E}
\newcommand{\PP}{\mathbb P}

\ifpdf
  \hypersetup{
    pdftitle={xxx}
  }
\fi

\hypersetup{colorlinks, allcolors = {MidnightBlue}}

\makeatletter
\renewcommand*{\defas}{%
\mathrel{\rlap{\raisebox{0.3ex}{$\m@th\cdot$}}\raisebox{-0.3ex}{$\m@th\cdot$}}=}
\makeatother

\makeatletter
\newtheorem*{rep@theorem}{\rep@title}
\newcommand{\newreptheorem}[2]{%
\newenvironment{rep#1}[1]{%
 \def\rep@title{#2 \ref*{##1}}%
 \begin{rep@theorem}}%
 {\end{rep@theorem}}}
\makeatother

\makeatletter
\DeclareRobustCommand\bigop[2][1]{%
  \mathop{\vphantom{\sum}\mathpalette\bigop@{{#1}{#2}}}\slimits@
}
\newcommand{\bigop@}[2]{\bigop@@#1#2}
\newcommand{\bigop@@}[3]{%
  \vcenter{%
    \sbox\z@{$#1\sum$}%
    \hbox{\resizebox{\ifx#1\displaystyle#2\fi\dimexpr\ht\z@+\dp\z@}{!}{$\m@th#3$}}%
  }%
}
\makeatother

\newtheorem{theorem}{Theorem}
\newtheorem{pro}[theorem]{Proposition}
\newreptheorem{pro}{Proposition}
\newtheorem{lem}[theorem]{Lemma}
\newreptheorem{lem}{Lemma}
\newtheorem{cor}[theorem]{Corollary}
\theoremstyle{definition}

\title{\vspace{-1cm}
{The asymptotic distribution of the $k$-Robinson-Foulds dissimilarity measure on labelled trees}}

\author [1] {Michael Fuchs}
\author[2]{Mike Steel}

\affil[1]{Department of Mathematical Sciences, National Chengchi University, Taipei 116, Taiwan}
\affil[2]{Biomathematics Research Centre, University of Canterbury, Christchurch, New~Zealand}

\begin{document}

\maketitle

\begin{abstract}
Motivated by applications in medical bioinformatics, Khayatian et al. (2024) introduced a family of metrics on Cayley trees (the $k$-RF distance, for $k=0, \ldots, n-2$)  and explored their distribution on pairs of random Cayley trees via simulations. In this paper, we investigate this distribution mathematically, and derive exact asymptotic descriptions  of the distribution of the $k$-RF metric for the extreme values $k=0$ and  $k=n-2$, as $n$ becomes large. We show that a linear transform of the $0$-RF metric converges to a Poisson distribution (with mean 2) whereas a similar transform for the $(n-2)$-RF metric leads to a normal distribution (with mean $\sim ne^{-2}$). These results (together with the case $k=1$ which behaves quite differently, and $k=n-3$) shed light on the earlier simulation results, and the predictions made concerning them. 
\end{abstract}

{\em Keywords:} Cayley trees,  $k$-RF dissimilarity metric, asymptotic enumeration, Poisson and Normal distributions.

\section{Introduction}

The Robinson-Foulds (RF) dissimilary measure is a simple and easily-computed metric on  phylogenetic trees (briefly, these are trees with labelled leaves and unlabelled interior vertices). Essentially, the RF metric measures how many `splits' (bipartitions of the labelled leaf set) are present in exactly one of the two trees (where a split results from deleting an edge in each tree). The RF metric is easily computed (unlike most other tree metrics in phylogenetics) and this has contributed to its wide  application in evolutionary biology; indeed, the 1981 paper that described it \cite{rob} has been cited more than 3000 times.  Moreover, the distribution of the RF metric on two random (binary or general) phylogenetic trees has been studied in the 1980s \cite{hen84, ste88}. 

In 2024, the authors of \cite{kvz24} investigated an analogous concept to define a class of metrics on the set of trees that have {\em all} their vertices labelled by distinct elements of some fixed set (i.e. the classical class of Cayley trees). This class of metrics,  which is parameterised by an integer $k\geq 0$ and denoted $k$-RF, is defined formally in the next section.   The motivation for this approach in \cite{kvz24} stems from its potential application in biomedical studies involving tumour cell evolution. 

A key part of the analysis in \cite{kvz24} was to estimate the distribution of the $k$-RF distance between two random (and independently sampled) Cayley trees. Their simulations led to various suggestions concerning possible limiting distributions as the number of vertices ($n$) grows. 

In this paper we investigate mathematically the distribution of the $k$-RF metric on Cayley trees, with a focus on asymptotic results, particularly for values of $k$ at (or near) their extreme values ($k=0$ and $k=n-2$). We then comment on how our results help interpret the findings suggested by simulations in \cite{kvz24}.

\subsection{Definitions}
A {\em Cayley tree on $[n]$} is a tree with vertex set $[n]=\{1,\dots, n\}$.  A classical result in combinatorics is that there are $n^{n-2}$ such trees \cite{cay}.  A {\em random Cayley tree} is a Cayley tree sampled uniformly at random from the set of all Cayley trees on $[n]$.  Note that a Cayley tree, as every tree, has $n-1$ edges.

The following measure was introduced and studied in the recent paper \cite{kvz24}.
Given a Cayley tree $T$, and $k \in \{0, 1, 2, \ldots, n-2\}$ a  {\em $k$-local split} of $T$ corresponding to edge $e = \{u,v\}$ is the pair of disjoint subsets $A,B$ of $[n]$ where $A$ (resp. $B$) is the subset of vertices of $T\setminus e$ within distance $k$ of $u$ (resp. within distance $k$ of $v$). Throughout the paper, we will use $A\vert B$ to denote such a split.

Notice in particular that the $0$-local split corresponding to edge $\{u, v\}$ is $\{u\}\vert \{v\}$. At the other extreme, an $(n-2)$-local split is the bipartition (i.e., partition into two blocks) of $[n]$ obtained by deleting an edge from the tree. 

\begin{figure}[htb]
\centering
\includegraphics[scale=0.5]{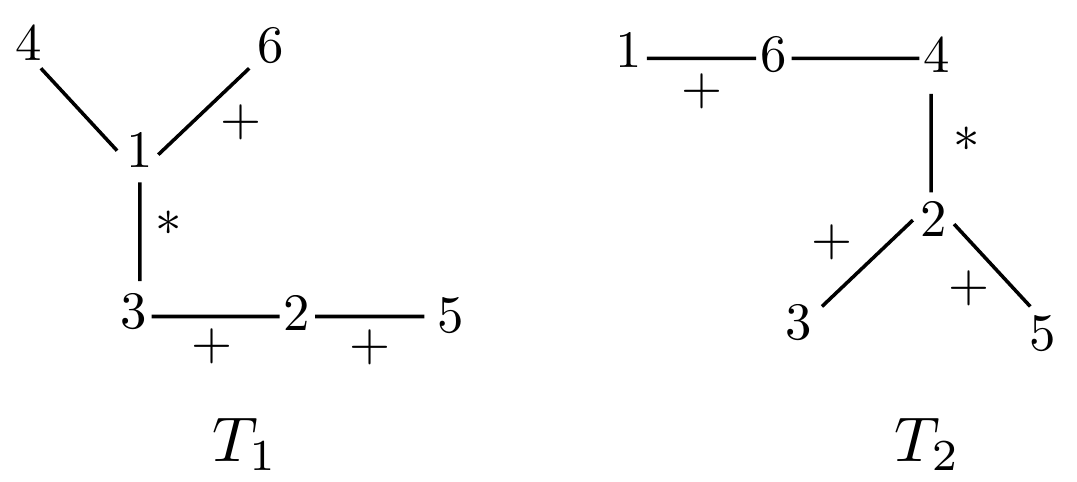}
\caption{The Cayley trees $T_1$ and $T_2$ on $[6]$ share three $0$-local splits (corresponding to the edges marked $+$) and one $3$-local split (corresponding to the edge marked $*$). Thus, $d_{0-RF}(T_1,T_2)= 4$ and $d_{(n-2)-RF}(T_1,T_2)= 8$.}
\label{fig1RF}
\end{figure}

Given two Cayley trees $T$ and $T'$ and $k \in \{0, \ldots, n-2\}$ the {\em $k$-RF distance} between $T$ and $T'$, denoted $d_{k-RF}(T, T')$, is the number of $k-$local splits they differ on.  Formally,
\begin{equation}\label{d-k-RF}
d_{k-RF}(T,T') = 2(n-1) - 2S_k(T, T'), 
\end{equation}
where $S_k(T,T')$ is the number of $k$-local splits shared by $T$ and $T'$. In other words, $d_{k-RF}(T,T')$ is the symmetric difference of the set of $k$-local splits of $T$ and the set of $k$-local splits of $T'$. Figure~\ref{fig1RF} provides an example with $k=0$ and $k=n-2$ to illustrate these definitions.

\vspace*{0.2cm}

Note that $d_{k-RF}$ is formally a metric on the set of Cayley trees on $[n]$, as established in Proposition 3 of \cite{kvz24} (the only metric axiom that requires a detailed argument is to establish that $d_{k-RF}(T,T')=0$ implies $T=T'$).
   
Note also that the restriction $k \leq n-2$ can be lifted, but there is no reason to do so, since any such metric would be identical to the case $k=n-2$ (since the maximal distance from any edge of a tree to any other vertex is $n-2$).

In this paper we first focus on the two  extreme cases for $k$ $(k=0, n-2$), and then consider the cases $k=1$ and $k=n-3$. 
We show that at one extreme ($k=n-2$) the number of  $k$-local splits shared by two random Cayley trees is asymptotically normal with a mean and variance that grows linearly with $n$, while at the other extreme ($k=0$) the number of $k$-local splits shared by two random Cayley trees converges to a Poisson distribution with mean $2$ as $n$ grows. The case $k=1$ exhibits a quite different behaviour, whereas the case $k=n-3$ is asymptotically the same as for $k=n-2.$

\section{The $(n-2)$--RF distribution}\label{n-2-RF}

Let $\T$ and $\T'$ be two independently sampled random Cayley trees on $[n]$,  let $S_n$ denote the number of $(n-2)$-local splits that $\T$ and $\T'$ share, and let $\tilde{S_n}$ denote the number of  $(n-2)$-local splits that $\T$ and $\T'$ share for which one side of the split is a single leaf. 

\bigskip

\begin{pro}\label{pro1} 
The following holds.
\begin{itemize}
    \item[(i)] $\EE(\tilde{S}_n(\tilde{S}_n-1) \cdots (\tilde{S}_n-k+1)) = k!\binom{n}{k} \left(1-\frac{k}{n}\right)^{2n-4}.$
    \item[(ii)] As $n\rightarrow \infty$,
    $\EE(\tilde{S}_n) \sim ne^{-2}$ and
    ${\rm Var}(\tilde{S}_n) \sim n(e^{-2}-3e^{-4})$
    and the same hold if $\tilde{S}$ is replaced by $S$.
\end{itemize}
\end{pro}
\begin{proof}
First observe that
$\EE(\tilde{S}_n(\tilde{S}_n-1)\cdots(\tilde{S}_n-k+1)) = k!\EE \left(\binom{\tilde{S}_n}{k}\right)$
and $\EE \left(\binom{\tilde{S}_n}{k}\right)$ is the sum over all $U \in \binom{[n]}{k}$ of the proportion of pairs of Cayley trees on $[n]$ for which  $U$ is a set of leaves in each tree.  For any set $U \in \binom{[n]}{k}$ this proportion is just:
$$\left(\frac{(n-k)^{n-k-2} \cdot (n-k)^k}{n^{n-2}}\right)^2 = (1-k/n)^{2n-4}.$$
Combining these observations gives Part (i).

The asymptotic expansion of the mean in part (ii) follows from part (i) by standard arguments. As for the variance, again part (ii) can be used since
\begin{equation}\label{var}
{\rm Var}(\tilde{S}_n)={\mathbb E}(\tilde{S}_n(\tilde{S}_n-1))+\EE(\tilde{S}_n)-(\EE(\tilde{S}_n))^2.
\end{equation}
However, the computation is slightly longer than that for mean; see the first proof of the next theorem. The proof of the final claim, namely, that the same expansions hold for $S_n$ follows from this proof, too.
\end{proof}

The main result of this section is that $S_n$ is asymptotically normally distributed as $n\rightarrow \infty$. In other words, if we normalise $S_n$ by subtracting its (asymptotic) mean and dividing by its (asymptotic) standard deviation, then the resulting random variable converges in distribution to the standard normal distribution $N(0,1)$, as follows.

\bigskip

\begin{theorem}\label{clt-Sn}
As $n \rightarrow \infty,$
\begin{equation}\label{clt-Sn-2}
    \frac{S_n - ne^{-2}}{\sqrt{(e^{-2}-3e^{-4})n}} \stackrel{d}{\longrightarrow} N(0,1).
\end{equation}
\end{theorem}

\begin{proof}

We first show that it suffices to establish this result with $S_n$ replaced by $\tilde{S}_n$. Observe that for $k \leq n/2$, the expected number of $(n-2)$-local splits that (i) partition $[n]$ into sets of size $k$ and $n-k$, and (ii) are shared by two independently  sampled random Cayley trees on $[n]$ is $$\binom{n}{k}\left( \frac{k^{k-1}(n-k)^{n-k-1}}{n^{n-2}}\right)^2.$$
For each fixed $k$, this expression is asymptotic  (as $n\rightarrow \infty$) to $ne^{-2}$ for $k=1$, to $2e^{-4}$ for $k=2$, and of order $n^{2-k}$ for $k>2$. This establishes our claim that it suffices to establish the result for $\tilde{S}_n$ in place of $S_n$,  since if $W_{n,k}$ denotes the number of $(n-2)$-local splits with both sides having  size at least $k$ shared by two independently sampled random Cayley trees on $[n]$, then $\PP(W_{n,k} \geq 1) \leq \EE(W_{n,k}).$

We now give two quite different proofs of Theorem~\ref{clt-Sn} using this last observation.  The first is combinatorial; it establishes  that the normalized central moments of $\tilde{S}_n$ converge to those of the standard normal distribution. The second argument is more explicitly probabilistic, relying on a known result concerning the asymptotic distribution of the number of leaves in a random Cayley tree.

 {\em Proof 1:} We first derive a (full) asymptotic expansion of the $k$-th factorial moment, where we use the closed-form expression from Proposition~\ref{pro1}(i). Note that 
 \begin{align}
 \left(1-\frac{k}{n}\right)^{2n-4}&=e^{(2n-4)\log(1-k/n)}=e^{-(2n-4)\sum_{\ell\geq 1}(k/n)^{\ell}}\nonumber\\
 &=e^{-2k-\sum_{\ell\geq 1}(2k-4)(k/n)^{\ell}}\nonumber\\
 &\sim e^{-2k}\sum_{\ell\geq 0}\frac{r_{\ell}(k)}{n^{\ell}},\label{asym-part-1}
 \end{align}
 where $r_{\ell}(k)$ is a polynomial of degree $2\ell$ with leading coefficient $(-1)^{\ell}/\ell!$ (and thus the asymptotic equivalence holds for $k=o(\sqrt{n}$)). Next,
 \begin{align}
 k!\binom{n}{k}&=n(n-1)\cdots(n-k+1)=\sum_{\ell\geq 0}s(k,\ell)n^{\ell}\nonumber\\
 &=n^k\sum_{\ell\geq 0}\frac{s(k,k-\ell)}{n^{\ell}},\label{asym-part-2}
 \end{align}
 where $s(k,k-\ell)$ are the Stirling numbers of first kind which are also polynomials of degree $2\ell$ and leading coefficient
 \[
 s(k,k-\ell)\sim\frac{(-1)^{\ell}k^{2\ell}}{2^{\ell}\ell!}.
 \]
Multiplying (\ref{asym-part-1}) and (\ref{asym-part-2}) gives
\begin{equation}\label{asym-fm}
{\mathbb E}(\tilde{S}_n(\tilde{S}_n-1)\cdots(\tilde{S}_n-k+1))\sim n^ke^{-2k}\sum_{\ell=0}^{\infty}\frac{p_{\ell}(k)}{n^{\ell}},
\end{equation}
where 
\[
p_{\ell}(k):=\sum_{j=0}^{\ell}s(k,k-j)r_{\ell-j}(k)
\]
are again polynomials in $k$ of degree $2\ell$. The leading coefficient of $p_{\ell}(k)$ is given by
\begin{equation}\label{lead-p}
\sum_{j=0}^{\ell}\frac{(-1)^{j}}{2^jj!}\cdot\frac{(-1)^{\ell-j}}{(\ell-j)!}=\frac{(-1)^{\ell}}{\ell!}\sum_{j=0}^{\ell}\binom{\ell}{j}\frac{1}{2^j}=\frac{(-1)^{\ell}}{\ell!}\left(\frac{3}{2}\right)^{\ell}.
\end{equation}

Our goal is to find the asymptotics of the $m$-th central moment of $\tilde{S}_n$, i.e.,
\begin{equation}\label{exp-cm}
{\mathbb E}(\tilde{S}_n-ne^{-2})^{m}=\sum_{j=0}^{m}\binom{m}{j}{\mathbb E}(\tilde{S}_n^j)(-1)^{m-j}n^{m-j}e^{-2(m-j)},
\end{equation}
where
\begin{equation}\label{exp-mom}
{\mathbb E}(\tilde{S}_n^{j})=\sum_{k=0}^{j}S(j,k){\mathbb E}(\tilde{S}_n(\tilde{S}_n-1)\cdots(\tilde{S_n}-k+1)).
\end{equation}
Here, $S(j,k)$ are the Stirling numbers of second kind. By plugging (\ref{asym-fm}) into (\ref{exp-mom}) and then plugging the resulting expression in turn into (\ref{exp-cm}), we obtain that
\[
{\mathbb E}(\tilde{S}_n-ne^{-2})^m\sim n^m\sum_{i=0}^{\infty}\frac{a_i}{n^{i}},
\]
where
\begin{equation}\label{ai}
a_i=\sum_{j=0}^{m}\binom{m}{j}(-1)^{m-j}e^{-2(m-j)}\sum_{k=j-i}^{j}S(j,k)e^{-2k}p_{k+i-j}(k).
\end{equation}

For fixed $i$, the number of terms of the last sum is also fixed. Pick one of the terms, say $k=j-i'$ with $0\leq i'\leq i$. Then, this term becomes
\[
S(j,j-i')e^{-2j+2i'}p_{i-i'}(j-i').
\]
As a function of $j$, $p_{i-i'}(j-i')$ is a polynomial in $j$ of degree $2i-2i'$ and $S(j,j-i')$ is a polynomial in $j$ of degree $2i'$. So, $q(j):=S(j,j-i')p_{i-i'}(j-i')$ is a polynomial in $j$ of degree $2i$. Overall, the contribution of $k=j-i'$ to $a_i$ becomes
\[
e^{-2(m-i')}\sum_{j=0}^{m}\binom{m}{j}(-1)^{m-j}q(j).
\]
This sum is $0$ whenever $i\leq(m-1)/2$ because of the following lemma.

\begin{lem}
For $s=0,1,\ldots$, we have
\[
\sum_{j=0}^m\binom{m}{j}(-1)^{m-j}j^{s}=\begin{cases} 0,&\text{if}\ s\leq m-1;\\ m!,&\text{if}\ s=m.\end{cases}
\]
\end{lem}
\begin{proof}
Consider
\[
\sum_{j=0}^{m}\binom{m}{j}(-1)^{m-j}z^j=(z-1)^{m}
\]
The claimed sum is obtained by applying $s$ times the operator $D:=z\frac{d}{dz}$ and then evaluating at $z=1$. Using induction, one easily shows that
\[
D^{s}(z-1)^m=(z-1)^{m-s}\lambda_s(z-1),
\]
where $\lambda_s(z)$ is a polynomial of degree $s$ with $\lambda_s(0)=m(m-1)\cdots(m-s+1)$. Evaluating at $z=1$ gives the claimed result.
\end{proof}

Consequently, from the arguments before the lemma, we have that $a_i=0$ for $i\leq (m-1)/2$ and thus
\begin{equation}\label{cm-O-estimate}
{\mathbb E}(\tilde{S}_n-ne^{-2})^{m}={\mathcal O}\left(n^{\lfloor m/2\rfloor}\right).
\end{equation}

We next refine this for $m=2m'$. Here, (\ref{ai}) with $i=m'$ becomes
\[
a_{m'}=\sum_{j=0}^{2m'}\binom{2m'}{j}(-1)^{2m'-j}e^{-2(2m'-j)}\sum_{k=j-m'}^{j}S(j,k)e^{-2k}p_{k+m'-j}(k).
\]
Set $k=j-m''$. Then, as explained above
\[
S(j,j-m'')p_{m'-m''}(j-m'')
\]
is a polynomial in $j$ of degree $2m'$ with leading coefficient
\[
\frac{1}{2^{m''}m''!}\cdot\frac{(-1)^{m'-m''}}{(m'-m'')!}\left(\frac{3}{2}\right)^{m'-m''}
\]
since $S(j,j-m'')$ is a polynomial in $j$ of degree $2m''$ with leading coefficient $1/(2^{m''}m''!)$ and because of (\ref{lead-p}). Thus, the contribution of $k=j-m''$ to $a_m'$ becomes
\[
e^{-2(2m'-m'')}\sum_{j=0}^{2m'}\binom{2m'}{j}(-1)^{2m'-j}S(j,j-m'')p_{m'-m''}(j-m'')
\]
which by the lemma equals
\[
e^{-2(2m'-m'')}(2m')!\frac{1}{2^{m''}m''!}\cdot\frac{(-1)^{m'-m''}}{(m'-m'')!}\left(\frac{3}{2}\right)^{m'-m''}.
\]
Summing this over $0\leq m''\leq m'$ gives
\begin{align*}
e^{-2m'}\frac{(2m')!}{2^{m'}}\sum_{m''=0}^{m'}\frac{(-3e^{-2})^{m'-m''}}{m''!(m'-m'')!}&=e^{-2m'}\frac{(2m')!}{2^{m'}m'!}\sum_{m''=0}^{m'}\binom{m'}{m''}(-3e^{-2})^{m'-m''}\\
&=e^{-2m'}(1-3e^{-2})^{m'}\frac{(2m')!}{2^{m'}m'!}.
\end{align*}
Overall,
\[
{\mathbb E}(\tilde{S}_n-ne^{-2})^{2m'}\sim e^{-2m'}(1-3e^{-2})^{m'}\frac{(2m')!}{2^{m'}m'!}n^{m'},
\]
or equivalently,
\[
{\mathbb E}\left(\frac{\tilde{S}_n-ne^{-2}}{\sqrt{(e^{-2}-3e^{-4})n}}\right)^{2m'}\longrightarrow\frac{(2m')!}{2^{m'}m'!}.
\]

Moreover, from (\ref{cm-O-estimate}),
\[
{\mathbb E}\left(\frac{\tilde{S}_n-ne^{-2}}{\sqrt{(e^{-2}-3e^{-4})n}}\right)^{2m'+1}\longrightarrow 0.
\]
The above two limits show that the moments of $(\tilde{S}_n-ne^{-2})/\sqrt{(e^{-2}-3e^{-4})n}$ converge to the moments of the standard normal distribution which concludes the proof of the claimed central limit theorem.

\bigskip

{\em Proof 2:} 
Recall that $\T$ and $\T'$ are two independently sampled random Cayley trees on $[n]$, and consider the resulting two  binary random variables $X_i, Y_i$ defined for each $i \in [n]$ as follows: Let 
$$X_i = 
\begin{cases}
    1, & \mbox{if $i$ is a leaf of $\T$;}\\
    0, & \mbox{otherwise}
\end{cases}
$$
and define $Y_i$ analogously, with $\T$ replaced by $\T'$. Then,
$\tilde{S}_n = \sum_{i=1}^nX_i\cdot Y_i$. 
We write this more conveniently as:
\begin{equation}
\label{sum}
    \tilde{S}_n = \sum_{i \in [n]:Y_i =1}X_i.
\end{equation}
Thus, $\tilde{S}_n$ is the sum of a random number $K_n = |\{i \in [n]:Y_i =1\}|$ of identically-distributed (but dependent) random variables. Notice, however, that the random variables $(X_i)$ are independent of the random variables $(Y_i)$ since the trees $\T$ and $\T'$ are sampled independently. Moreover, the $X_i$ (and also the $Y_i$) are exchangable random variables.

We now apply the following result which is a consequence of Theorem 2.3 of \cite{drm}: The number of leaves in a uniformly-sampled Cayley tree on $[n]$ is asymptotically normally distributed with  mean  and variance asymptotic to $e^{-1}n$ and $(e^{-2}-2e^{-4})n$, respectively.  Thus, $K_n/n$ converges almost surely to $e^{-1}$, and so   $\tilde{S}_n$ is asymptotically normal. The asymptotic mean and variance of $\tilde{S}_n$ is provided by Proposition~\ref{pro1}(ii).\end{proof}

{\bf Remark:} The result used in the second proof (the consequence of Theorem 2.3 of \cite{drm}) can be also proved with our method used in the first proof.

\bigskip

\begin{cor}
If $\T$ and $\T'$ are two independent Cayley trees on $[n]$ then $$\frac{d_{(n-2)-RF}(\T, \T')- 2n(1-e^{-2})}{2\sqrt{(e^{-2}-3e^{-4})n}}$$ converges in distribution to $N(0,1)$ as $n$ grows.
\end{cor}

\bigskip

\begin{proof}
By (\ref{d-k-RF}), 
\[
S_n=n-1-\frac{d_{(n-2)-RF}(\T,\T')}{2}.
\]
Plugging this into (\ref{clt-Sn-2}) and multiplying by $-1$ gives the desired result.
\end{proof}

\section{The $0$--RF distribution}\label{0-RF}

In this section, we show that the number of $0$-local splits shared by two independent Cayley trees on $[n]$ converges to a Poisson distribution with mean 2 as $n$ becomes large. Recall that a shared $0$-local split is a pair of elements of $[n]$ that correspond to an edge in both trees.

Our approach relies on the following result (Theorem 6.1 of \cite{moo}; an alternative proof of this theorem was described recently  in \cite{cam}).

\bigskip

\begin{pro}
\label{pro-span}
      The number of spanning trees of the complete graph on $[n]$ that contain a given spanning forest $F$ of $m$ trees is given by $$q_1 \cdots q_m \cdot n^{m-2},$$ where $q_i$ is the number of vertices in tree $i$.
\end{pro}

By setting  $q_i=2$ for $k$ values of $i$ and $q_i=1$ otherwise in this last proposition we obtain the following.

\bigskip

\begin{cor}\label{cor-span}
Suppose that $C$ is a collection of $j$ 2-sets of $[n]$ (i.e. sets of size 2). 
If 2-sets in $C$ are pairwise disjoint then the proportion of Cayley trees on $[n]$ that contain these 2-sets as edges is exactly 
$(\frac{2}{n})^j$. 
\end{cor}

Next, let  $S'_n$ be the number of pairs ($\{x,y\} \in \binom{[n]}{2}$) shared by  two independently sampled random Cayley trees on $[n]$.

\bigskip

\begin{cor}
\label{cor-span2}
$\EE(S_n') = 2(n-1)/n \sim 2$ and ${\rm Var}(S'_n) \sim 2.$
\end{cor}

\begin{proof}
For the mean, by Corollary~\ref{cor-span} (with $j=1$), we obtain
\[
\EE(S'_n) = \binom{n}{2}\left(\frac{2}{n}\right)^2 = \frac{2(n-1)}{n} \sim 2
\]
which shows the claim. 

For the variance, again by Corollary~\ref{cor-span} (with $j=2$), 
\[
\EE(S'_n(S'_{n}-1)) \sim\binom{n}{2}\binom{n-2}{2}\left(\frac{2}{n}\right)^4 \sim 4. 
\]
Combining this with formula (\ref{var}) (with $\tilde{S}_n$ replaced by $S_n'$) and $\EE(S'_n) \sim 2$ gives the claim also in this case.
\end{proof}

{\bf Remark:} 
It is easily shown that the expected number of edges shared by a random Cayley tree on $[n]$ and a  fixed Cayley tree on $[n]$ is also $2(n-1)/n.$ This invariance of the mean under fixing one of the two trees for $0$--RF  does not hold for $(n-2)$--RF, since trees with many leaves will tend to share more $(n-2)$-local splits with a random tree. Nevertheless, the distribution of the number of edges shared by a random Cayley tree on $[n]$ and a fixed Cayley tree $T_n$ on $[n]$ can be quite different than that arising from two random trees.  To see this, let $T_n$ have exactly one non-leaf vertex. Then every Cayley tree on $[n]$ shares at least one edge with $T_n$. Thus, the probability that a random Cayley tree on $[n]$ shares no edge with $T_n$ equals $0$ (for all $n$), and so is not described by a Poisson distribution, which we will show applies asymptotically for two random Cayley trees.

The main result of this section is the following.

\bigskip

\begin{theorem}
    The total variational distance between $S'_n$ and a Poisson random variable with mean $2$ converges to 0 (at rate $O(1/n)$) as $n$ grows. Thus, if $\T$ and $\T'$ are two independent random Cayley trees on $[n]$, then $d_{0-RF}(\mathcal{T}, \mathcal{T}')/2-(n-1)$ converges in distribution to a Poisson random variable with mean $2$.
\end{theorem}
\label{prostein}
    
\begin{proof}
  For $\{u,v\}$ in $\binom{[n]}{2}$, let $N_{uv}$ be the collection of sets $\{w,z\} \in \binom{[n]}{2}$ for which $|\{u,v\} \cap \{w,z\}|=1$. Clearly,
  \begin{equation}
      \label{counteq}
      |N_{uv}|= 2n-3.
  \end{equation}
Let $X_{uv}$ be the Bernoulli random variable that takes the value $1$ if $uv$ is an edge in two Cayley trees sampled uniformly at random, and is zero otherwise, and let  $p_{uv} = \PP(X_{uv}=1) = \EE(X_{uv})$. By Corollary~\ref{cor-span},

\begin{equation}
\label{ouv}
    p_{uv} = (2/n)^2
\end{equation}

{\bf Claim:}  For each $\{u,v\} \in \binom{[n]}{2}$, the random variable $X_{uv}$ is independent of the set of variables $(X_{rs}: \{r,s\} \in \binom{[n]}{2} \setminus N_{uv})$.

 {\em Proof of claim:} It suffices to show that 
 \begin{equation}
 \label{indep}
     \PP(\{X_{uv}=1\}\cap {\mathcal A}) = \PP(X_{uv}=1)\cdot\PP( {\mathcal A}) 
 \end{equation}
 for any assignment $\mathcal{A}$ of a state ($0$ or $1$) to each variable $X_{rs}$ in $\binom{[n]}{2} \setminus N_{uv}$.    Note that any such assignment corresponds to a graph  on $[n]\setminus \{u,v\}$, which, moreover, contains no cycles (due to how the random variables $X_*$ are defined) and so this graph is a forest. If this forest consists of trees of sizes $q_1, \ldots q_r$ then by  Proposition~\ref{pro-span}, the left-hand side of (\ref{indep}) is:
 \begin{equation}
 \left(\frac{2q_1\cdots q_r n^{r-1}}{n^{n-2}}\right)^2    
 \end{equation}
and the right-hand side of (\ref{indep}) is given by:
  \begin{equation}
\left(\frac{2}{n}\right)^2 \cdot \left(\frac{1 \cdot 1 \cdot q_1 \cdots q_r n^{r}}{n^{n-2}}\right)^2    
 \end{equation}
and these are identical expressions, and so (\ref{indep}) holds.   This completes the proof of the Claim.

Next, let $\lambda_n = \EE(S'_n) = 2(1-\frac{1}{n})$, and let $Po(\lambda_n)$ denote a Poisson random variable with mean $\lambda_n$.  We now apply a version of the Stein-Chen bound (the `dissociated case' (see e.g. Theorem 4.4.8 of \cite{roc})) which states  that:
\begin{equation}
\label{chen}
    d_{TV}(S'_n, Po(\lambda_n)) \leq \min\left(1, \frac{1}{\lambda_n}\right)\sum_{uv \in \binom{[n]}{2}}\left(p_{uv}^2 + \sum_{u'v' \in N_{uv}} (p_{uv}p_{u'v'} +\mathbb{E}(X_{uv}X_{u'v'}))\right).
\end{equation}

On the right of this inequality, the terms  $p_{uv}^2$  and $p_{uv}p_{u'v'}$ equal $(2/n)^4$ by (\ref{ouv}).
The term $\EE(X_{uv}X_{u'v'})= \PP(X_{uv} =1 \wedge X_{u'v'}=1)$ is equal to $(3/n^2)^2$  since $\{u,v\}$ and $\{u',v'\}$ have exactly one vertex in common and so the proportion of Cayley trees on $[n]$ that share this 3-vertex tree is $3/n^2$ by Proposition~\ref{pro-span}.
By (\ref{counteq}) it follows that the second summation terms on the right of Inequality (\ref{chen})  are of order $n^{-3}$ and since the outer sum is over $\binom{n}{2}$ terms this provides a bound of order $n^{-1}$ as claimed. 

Theorem~\ref{prostein} now follows by noting that (i) the variational distance between a Poisson random variable with mean $\lambda_n$ and one with mean 2 tends to 0 at rate $n^{-1}$ as $n$ grows, and (ii) application of the triangle inequality to this metric.
\end{proof}

{\bf Remark:} The following weaker result can be  established by more elementary means:  

\bigskip

\begin{pro}
\label{pro2}
    For each integer $r\geq 0$, the probability that $S'_n=r$ converges to the Poisson probability $e^{-2} 2^r/r!$ as $n \rightarrow \infty.$ 
\end{pro}

\vspace*{0.2cm}

To establish this result,  we first require the following lemma.

\bigskip

\begin{lem}
\label{lem-span}
   If $k$ 2-sets are selected uniformly at random from $\binom{[n]}{2}$ the probability that these $2$-sets are pairwise disjoint tends to 1 as $n \rightarrow \infty$. 
\end{lem}
\begin{proof}
There are $$\frac{1}{k!} \times \binom{n}{2} \times \binom{n-2}{2} \times \cdots \times \binom{n-2k+2}{2}$$
ways to select $k$ disjoint pairs of elements from $[n]$, and $$\frac{1}{k!} \times \binom{n}{2} \times \left(\binom{n}{2}-1\right) \times \cdots \times \left(\binom{n}{2}-k+1\right)$$ ways to select $k$ distinct pairs of elements from $[n]$, and the ratio of the first product to the second product is $1-o(1)$ as $n\rightarrow \infty$.
\end{proof}

\begin{proof}[Proof of Proposition~\ref{pro2}.]
We apply the (generalised) Principle of Inclusion and Exclusion, in conjunction with Lemma~\ref{lem-span} and Corollary~\ref{cor-span}. This implies that for $r=0, 1, 2, \ldots$, we have:
$${\mathbb P} (S'_n=r) =\sum_{m=r}^{n-1}(-1)^{m-r}\binom{m}{r}\frac{1}{m!}\binom{n}{2}\binom{n-2}{2}\cdots\binom{n-2m+2}{2}\left(\frac{2}{n}\right)^{2m} + o(1),$$ 
where $o(1)$ is a term that tends to $0$ as $n$ grows (by Lemma~\ref{lem-span}).
The summation term  is asymptotic to $$\frac{1}{r!}\sum_{m=r}^{n-1} (-1)^{m-r}\frac{2^m}{(m-r)!}$$ 
as $n \rightarrow \infty$. Thus, for each $r$, ${\mathbb P} (S'_n=r) \rightarrow e^{-2}2^r/r!$ as $n \rightarrow \infty.$
\end{proof}

\subsection*{A related result}
Suppose we select a set $P_n$ of $n-1$ pairs of elements of $[n]$ uniformly at random (and without replacement).  Such a set corresponds to a graph on $[n]$ but in general this graph will not correspond to the set of edges of a Cayley tree, since the graph may be disconnected and/or contain cycles. Nevertheless,  for a fixed (or random) Cayley tree $T_n$ with $n$ leaves, let $S^*_{T_n}$ be the number of edges of $T_n$ that are present in $P_n$. Then we have the following result.

\bigskip

\begin{pro}
   As $n \rightarrow \infty$, $S^*_{T_n}$ converges in distribution to a Poisson random variable with mean $2$.
\end{pro}

\begin{proof}
We use the following result from \cite{mst} concerning a Poisson limit of the hypergeometric distribution. For $k= 0,1,2, \ldots,$ we have:
$$\lim_{r,s,t \rightarrow \infty; \frac{r\cdot t}{r+s} \rightarrow \lambda} \frac{\binom{r}{k}\binom{s}{t-k}}{\binom{r+s}{t}} = \frac{\lambda^ke^{-\lambda}}{k!}$$
with $r=t=n-1$ and $s=\binom{n}{2} -r$, which satisfy the conditions under the limit sign, with $\lambda = 2.$ Note that we are sampling $n-1$ pairs of elements of $[n]$  from the set of  all $r+s = \binom{n}{2}$ pairs, and so the probability on the left is the probability of sampling (without replacement) exactly $k$ of the edges in $T_n$.
\end{proof}

\section{The $(n-3)$--RF and $1$--RF distributions}

\paragraph{The $(n-3)$--RF distribution.} We first observe that the $(n-2)$-local splits and $(n-3)$-local splits are almost the same: the former consist of any partition of the set $[n]$ into two blocks whereas the latter have in addition the splits $A\vert B$ with $A$ a singleton and $B$ consisting of $[n]\setminus A$ with one more element removed.

Let now $\T$ and $\T'$ be as in Section~\ref{n-2-RF} and $U_n$ resp $\tilde{U}_n$ be similar to $S_n$ and $\tilde{S}_n$ also from Section~\ref{n-2-RF} only with $(n-2)$-local splits in their definitions replaced by $(n-3)$-local splits. 

We first show the following lemma.

\bigskip

\begin{lem}\label{diff-to-0}
As $n\rightarrow\infty$,
\[
{\mathbb E}(\vert S_n-U_n\vert)=o(1).
\]
\end{lem}
\begin{proof}
First, observe that 
\[
{\mathbb E}(\vert S_n-U_n\vert)={\mathbb E}(\vert\tilde{S}_n-\tilde{U}_n\vert)
\]
and thus we can concentrate on the splits $A\vert B$ where $A$ is a singleton. 

As mentioned above, these splits are either of the form $\{x\}\vert [n]\setminus\{x\}$ or of the form $\{x\}\vert [n]\setminus\{x,y\}$ where $x\ne y$. The probabilities that $\T$ and $\T'$ contain the former are $4/n^2$ and $(2/n-n!/n^{n-1})^2$ when considering $(n-2)$-local splits and $(n-3)$-local splits, respectively. On the other hand, the latter have just to be considered in the case of $(n-3)$-splits where they are contained in $\T$ and $\T'$ with probability $(n-2)!^2/n^{2n-4}$. Thus, we have
\[
{\mathbb E}(\vert\tilde{S}_n-\tilde{U_n}\vert)\leq n\left(\frac{4}{n^2}-\left(\frac{2}{n}-\frac{n!}{n^{n-1}}\right)^2\right)+\frac{n(n-2)!^2}{n^{2n-4}}\sim\frac{4\sqrt{2\pi n^3}}{e^n},
\]
where we used Stirling's formula in the last step. This implies the claim.
\end{proof}
From this lemma, we obtain the following theorem.

\bigskip

\begin{theorem}
    As $n \rightarrow \infty,$
    $$\frac{U_n - ne^{-2}}{\sqrt{(e^{-2}-3e^{-4})n}} \stackrel{d}{\longrightarrow} N(0,1)$$
and consequently, for two independent Cayley trees $\T$ and $\T'$ on [n],
\[
\frac{d_{(n-3)-RF}(\T, \T')- 2n(1-e^{-2})}{2\sqrt{(e^{-2}-3e^{-4})n}}\stackrel{d}{\longrightarrow}N(0,1).
\]
\end{theorem}
\begin{proof}
We have,
\[
\frac{U_n-ne^{-2}}{\sqrt{(e^{-2}-3e^{-4})n}}=\frac{U_n-S_n}{\sqrt{(e^{-2}-3e^{-4})n}}+\frac{S_n-ne^{-2}}{\sqrt{(e^{-2}-3e^{-4})n}}.
\]
From Lemma~\ref{diff-to-0} and Markov's inequality, we obtain the following convergence in probability:
\[
\frac{U_n-S_n}{\sqrt{(e^{-2}-3e^{-4})n}}\stackrel{P}{\longrightarrow} 0.
\]

Thus, the claim follows from Theorem~\ref{clt-Sn} and Slutsky's theorem.
\end{proof}

\paragraph{The 1-RF distribution.} The $1$-local splits $A\vert B$ consists of all disjoint subsets of the set $[n]$ with either $\vert A\vert$ is a singleton and 
$\vert B\vert\geq 2$ or $\vert A\vert,\vert B\vert\geq 2$. Let the former be called of {\it type-1} and the latter of {\it type-2}. We intend to count the number of Cayley trees which contain a fixed such split. We start with a proposition.

\bigskip

\begin{pro}\label{count-forests}
The number of labeled, ordered forests with $n$ vertices consisting of $s\leq n$ trees where each tree contains label $1$ and all other labels are different and from the set $\{2,\ldots,n-s\}$ equals $sn^{n-s-1}$. 
\end{pro}
\begin{proof}
Since the number of Cayley trees on $[n]$ equals $n^{n-2}$, the exponential generating function of the number of Cayley trees on $[n]$ where the vertex of label $1$ does not contribute to its size is given by
\[
\sum_{n\geq 1}n^{n-2}\frac{z^{n-1}}{(n-1)!}=:\frac{A(z)}{z}.
\]
Here, $A(z)$ is the exponential generating function of the number of rooted Cayley trees on $[n]$ (which are counted by $n^{n-1}$). Recall that
\[
A(z)=ze^{A(z)},
\]
which follows from the exponential formula; see, e.g., Chapter~14 in \cite{LiWi}. The lemma asks now for the following count:
\begin{align*}
(n-s)![z^{n-s}]\left(\frac{A(z)}{z}\right)^s&=(n-s)![z^{n}](A(z))^s\\
&=\frac{s(n-s)!}{n}[\omega^{n-s}]e^{n\omega}=sn^{n-s-1},
\end{align*}
where we used Lagrange inversion in the second step. This proves the claimed result.
\end{proof}

\begin{figure}
\centering
\begin{tikzpicture}
\filldraw[black] (0cm,0cm) circle (1.8pt) node[label=below:{$\alpha$}] {};
\filldraw[black] (1.2cm,0cm) circle (1.8pt) node[label=below:{$\beta$}] {};
\filldraw[black] (2.4cm,0.8cm) circle (1.8pt) node[label=right:{$x_1$}] {};
\filldraw[black] (2.4cm,0.4cm) circle (1.8pt) node[label=right:{$x_2$}] {};
\filldraw[black] (2.4cm,-0.4cm) circle (1.8pt) node[label=right:{$x_{k-1}$}] {};
\filldraw[black] (2.4cm,-0.8cm) circle (1.8pt) node[label=right:{$x_k$}] {};

\draw[-,line width=0.7pt] (0cm,0cm) -- (1.2cm,0cm);
\draw[-,line width=0.7pt] (1.2cm,0cm) -- (2.4cm,-0.8cm);
\draw[-,line width=0.7pt] (1.2cm,0cm) -- (2.4cm,-0.4cm);
\draw[dotted,line width=0.7pt] (2.4cm,-0.15cm) -- (2.4cm,0.15cm);
\draw[-,line width=0.7pt] (1.2cm,0cm) -- (2.4cm,0.4cm);
\draw[-,line width=0.7pt] (1.2cm,0cm) -- (2.4cm,0.8cm);

\filldraw[black] (6.8cm,0.8cm) circle (1.8pt) node[label=left:{$y_1$}] {};
\filldraw[black] (6.8cm,0.4cm) circle (1.8pt) node[label=left:{$y_2$}] {};
\filldraw[black] (6.8cm,-0.4cm) circle (1.8pt) node[label=left:{$y_{\ell-1}$}] {};
\filldraw[black] (6.8cm,-0.8cm) circle (1.8pt) node[label=left:{$y_{\ell}$}] {};
\filldraw[black] (8cm,0cm) circle (1.8pt) node[label=below:{$\alpha$}] {};
\filldraw[black] (9.2cm,0cm) circle (1.8pt) node[label=below:{$\beta$}] {};
\filldraw[black] (10.4cm,0.8cm) circle (1.8pt) node[label=right:{$x_1$}] {};
\filldraw[black] (10.4cm,0.4cm) circle (1.8pt) node[label=right:{$x_2$}] {};
\filldraw[black] (10.4cm,-0.4cm) circle (1.8pt) node[label=right:{$x_{k-1}$}] {};
\filldraw[black] (10.4cm,-0.8cm) circle (1.8pt) node[label=right:{$x_k$}] {};

\draw[-,line width=0.7pt] (6.8cm,-0.8cm) -- (8cm,0cm);
\draw[-,line width=0.7pt] (6.8cm,-0.4cm) -- (8cm,0cm);
\draw[dotted,line width=0.7pt] (6.8cm,-0.15cm) -- (6.8cm,0.15cm);
\draw[-,line width=0.7pt] (6.8cm,0.4cm) -- (8cm,0cm);
\draw[-,line width=0.7pt] (6.8cm,0.8cm) -- (8cm,0cm);
\draw[-,line width=0.7pt] (8cm,0cm) -- (9.2cm,0cm);
\draw[-,line width=0.7pt] (9.2cm,0cm) -- (10.4cm,-0.8cm);
\draw[-,line width=0.7pt] (9.2cm,0cm) -- (10.4cm,-0.4cm);
\draw[dotted,line width=0.7pt] (10.4cm,-0.15cm) -- (10.4cm,0.15cm);
\draw[-,line width=0.7pt] (9.2cm,0cm) -- (10.4cm,0.4cm);
\draw[-,line width=0.7pt] (9.2cm,0cm) -- (10.4cm,0.8cm);

\draw (1.2cm,-1.3cm) node {(i)};
\draw (8.6cm,-1.3cm) node {(ii)};
\end{tikzpicture}
\caption{Type-1 (left) and Type-2 (right) $1$-local splits.}\label{split-types}
\end{figure}
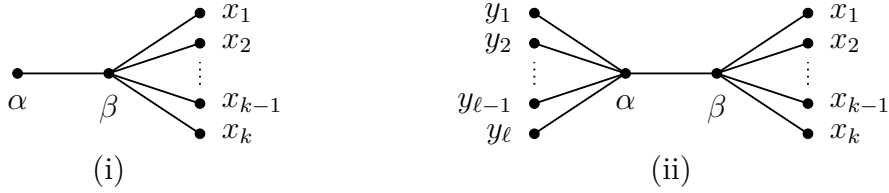

\begin{cor}\label{count-splits}
\begin{itemize}
\item[(i)] The number of Cayley trees which contain a fixed type-1 split $A\vert B$ with $\vert B\vert=k+1$ equals $k(n-2)^{n-k-3}$.
\item[(ii)] The number of Cayely trees which contain a fixed type-2 split $A\vert B$ with $\vert A\vert=~\ell+1$ and $\vert B\vert=k+1$ equals $(\ell+k)(n-2)^{n-\ell-k-3}$.
\end{itemize}
\end{cor}
\begin{proof}
Fix a type-1 split; see the left part of Figure~\ref{split-types}. Next, fix a labeled, ordered forest with $n-2$ vertices consisting of $k$ trees, each of them having the label $1$ and having all their other labels different. Attach them (in order) to $x_1,\ldots,x_k$ at the vertices of label $1$ which then gets replaced by $x_1,\ldots,x_k$. Finally, relabel all vertices except $\alpha,\beta$, and $x_1,\ldots,x_k$ in an order-consistent way such that the final tree has labels $1,\ldots,n$ (i.e., is a Cayley tree on $[n]$). Since there are $k(n-2)^{n-k-3}$ possible choices for the forest in the above procedure (Lemma~\ref{count-forests}) and all Cayley trees containing the fixed split are (uniquely) obtained in this way, we are done.

The proof for a type-2 split works similarly.
\end{proof}

Now let $\T$ and $\T'$ be two independently sampled random Cayley trees on $[n]$. Similarly to Section~\ref{0-RF}, denote by $U_n'$ the splits shared by $\T$ and $\T'$. Then, we have the following result.

\bigskip

\begin{theorem}
We have,
\[
{\mathbb E}(U_n')={\mathcal O}(1/n).
\]
Thus, the limiting distribution of $U_n'$ is degenerate at $0$.
\end{theorem}
\begin{proof}
Let $X_S$ the indicator random variable which counts whether or not $\T$ contains split $S$ and define $Y_S$ similarly for $\T'$. Then,
\[
U_n'=\sum_{S}X_S\cdot Y_S,
\]
where the sum runs over all $1$-local splits $S$. Next, denote by $p_S$, the probability that $X_S=1$. Then, by Corollary~\ref{count-splits},
\[
p_S=\begin{cases}k(n-2)^{n-k-3}n^{-n+2}\leq kn^{-k-1},&\text{if $S$ is as in Figure~\ref{split-types},(i)};\\
(\ell+k)(n-2)^{n-\ell-k-3}n^{-n+2}\leq (\ell+k)n^{-\ell-k-1},&\text{if $S$ is as in Figure~\ref{split-types},(ii).}\end{cases}
\]
Note that ${\mathbb E}(X_SY_S)=p_{S}^2$. Thus, in order to compute ${\mathbb E}(U_n')$, we need to know the number of splits $S$ of type-1 and type-2 as depicted in Figure~\ref{split-types}. Let $N_1(k)$ be the number of those on the left of the figure and $N_2(\ell,k)$ the number of those on the right. Then,
\[
N_1(k)=n(n-1)\binom{n-2}{k}\leq\frac{n^{k+2}}{k!}
\]
as there are $n(n-1)$ choices for $\alpha$ and $\beta$ and $\binom{n-2}{k}$ choices for $x_1,\ldots,x_k$. Likewise
\[
N_2(\ell,k)=\frac{1}{2}\binom{n}{2}\binom{n-2}{\ell}\binom{n-2-\ell}{k}\leq\frac{n^{\ell+k+2}}{4\ell!k!}
\]
Thus, we obtain for the expected value of $U_n'$
\begin{align*}
{\mathbb E}(U_n')&\leq \sum_{k\geq 1}N_1(k)k^2n^{-2k-2}+\sum_{\ell,k\geq 1}N_2(\ell,k)(l+k)^2n^{-2\ell-2k-2}\\
&\leq\sum_{k\geq 1}\frac{k}{(k-1)!n^{k}}+\sum_{\ell,k\geq 1}\frac{(\ell+k)^2}{4\ell!k!n^{\ell+k}}={\mathcal O}\left(\frac{1}{n}\right),
\end{align*}
where we used the above estimates of the probability $p_S$. This proves the claimed result.
\end{proof}

As a corollary of the last theorem, we obtain the following.

\bigskip

\begin{cor}
For two independently sampled random Cayley trees $\T$ and $\T'$ on $[n]$, the limiting distribution of $d_{1-RF}(\T,\T')$ is degenerate at $2n-2$.
\end{cor}

\section{Concluding comments}
Our theoretical results shed light on the shape of the empirical distribution of $k$-RF for $n=4,5,6,7$ in \cite{kvz24} obtained by simulations and presented
 in the histograms of Figs. 3--6 of that paper. Firstly, our theoretical results are consistent with the general shape of these empirical distributions. Second,  the   authors of  \cite{kvz24} state (in Section 4.2) that ``it seems the pairwise $0$-RF and $(n-2)$-RF scores have a Poisson distribution''.  While a Poisson distribution provides a moderate fit for their empirical study based on small values of $n$, our results show that as $n$ increases, the distributions become somewhat different. Specifically, the $0$-RF distance does not converge to a Poisson distribution, however, $(d_{0-RF}-2(n-1))/2$ does converge to a Poisson distribution (with mean $2$). On the other hand, $d_{(n-2)-RF}$ converges to a normal distribution with a mean that is different from its variance, and so $d_{(n-2)-RF}$ is not described by a Poisson distribution when $n$ is large.  
 
Finally, our analysis of the $k$-RF metric for $k=1$ explains why the histograms for this value of $k$ in \cite{kvz24} appear so different to the cases $k=0$ and $k=n-2$. By contrast, the case $k=n-3$ is asymptotically the same as for $k=n-2$ even though this is not (yet) visible in the histograms from \cite{kvz24}.  It may also be of interest  to investigate the distribution of the $k$-RF metric for other values of $k$.

\section{Acknowledgments} This research was carried out during a sabbatical stay of MF at the Biomathematics Research Center, University of Canterbury, Christchurch. He thanks the department and MS for hospitality, and the National Science and Technology Council, Taiwan (research grants NSTC-113-2918-I-004-001 and NSTC-113-2115-M-004-004-MY3) for financial support. MS thanks the NZ Marsden Fund for research support (23-UOC-003).

\end{document}